\documentclass[11pt]{amsart}
\usepackage{amssymb}
\usepackage{amscd}
\usepackage[all]{xy}
\xyoption{ps}
\xyoption{dvips}
\usepackage{amsthm}

\newtheorem{Thm}{Theorem}[section]

\newtheorem{Cor}[Thm]{Corollary}
\newtheorem{Prop}[Thm]{Proposition}
\newtheorem{Conj}[Thm]{Conjecture}
\newtheorem{Qu}[Thm]{Question}
\newtheorem{Prob}[Thm]{Problem}

\newcommand{\A}{\mathbb{A}}
\newcommand{\D}{\mathbb{D}}
\newcommand{\E}{\mathbb{E}}

\newcommand{\Z}{\mathbb{Z}}
\newcommand{\Q}{\mathbb{Q}}
\newcommand{\N}{\mathbb{N}}

\newcommand{\md}{\operatorname{mod}}

\newcommand{\Ext}{\operatorname{Ext}}

\newcommand{\M}{{\rm M}}
\newcommand{\GL}{{\rm GL}}
\newcommand{\field}{k}
\newcommand{\Or}{\mathcal O}
\newcommand{\rep}{\operatorname{rep}}

\begin{document}

\title{Extension-orthogonal components of nilpotent varieties}
\author{Christof Gei{\ss}}
\address{Christof Gei{\ss}\newline
Instituto de Matem\'aticas, UNAM\newline
Ciudad Universitaria\newline
04510 Mexico D.F.\newline
Mexico}
\email{christof@matem.unam.mx}

\author{Jan Schr\"oer}
\address{Jan Schr\"oer\newline
Department of Pure Mathematics\newline
University of Leeds\newline
Leeds LS2 9JT\newline
UK}
\email{jschroer@maths.leeds.ac.uk}


\begin{abstract}
Let $Q$ be a 
Dynkin quiver, and let $\Lambda$ be the corresponding 
preprojective algebra.
Let ${\mathcal I} = \{ C_i \mid i \in I \}$ be a set of pairwise different
indecomposable irreducible components of varieties of
$\Lambda$-modules such that generically there are no extensions
between $C_i$ and $C_j$ for all $i,j$.
We show that the number of elements in ${\mathcal I}$ is at most
the number of positive roots of $Q$. 
Furthermore, we give a module theoretic interpretation of
Leclerc's counterexample to a conjecture of Berenstein and
Zelevinsky.
\end{abstract}

\maketitle


\section{Introduction}


Let $\field$ be an algebraically closed field.
For a finitely generated $\field$-algebra $A$ let
$\md_A({\bf d})$ be the affine variety of $A$-modules with dimension
vector ${\bf d}$.
For irreducible components
$C_1 \subseteq \md_A({\bf d}_1)$ and
$C_2 \subseteq \md_A({\bf d}_2)$ 
define
\[
{\rm ext}^1_A(C_1,C_2) = \min \{ {\rm dim}\, \Ext^1_A(M_1,M_2) 
\mid (M_1,M_2) \in C_1 \times C_2 \}.
\]
An irreducible component $C \subseteq \md_A({\bf d})$ is {\it indecomposable}
if it contains a dense subset of indecomposable $A$-modules.
A general theory of irreducible components and their decomposition
into indecomposable irreducible components was developed in \cite{CBSc}.
Our aim is to apply this to Lusztig's nilpotent varieties.

If not mentioned otherwise, we always assume that $Q$ is a Dynkin quiver of
type $\A_n$, $\D_n$ or $\E_{6,7,8}$.
By  $R^+$ we denote the set of positive roots of $Q$, and by 
$\Lambda$ we denote the preprojective algebra associated to $Q$,
see \cite{Ri2}.
Let $n$ be the number of vertices of $Q$, and let
$\Lambda({\bf d}) = \md_{\Lambda}({\bf d})$, ${\bf d} \in \N^n$,
be the variety of $\Lambda$-modules with dimension vector ${\bf d}$.
The varieties $\Lambda({\bf d})$ are called 
{\it nilpotent varieties}.
We refer to \cite[Section 12]{L2} for basic properties.
Throughout, we only consider finite-dimensional modules.
Our main result is the following:

\vspace{0.2cm}\begin{Thm}\label{mainresult}
Assume that $\{ C_i \subseteq \Lambda({\bf d}_i) \mid i \in I \}$
is a set of pairwise different indecomposable irreducible components 
such that
${\rm ext}_{\Lambda}^1(C_i,C_j) = 0$ for all $i,j \in I$.
Then $|I| \leq |R^+|$.
\end{Thm}\vspace{0.2cm}

As a consequence we get the following result on $\Lambda$-modules without
self-extensions:

\vspace{0.2cm}\begin{Cor}\label{maincor}
Let $M$ be a $\Lambda$-module
with ${\rm Ext}_{\Lambda}^1(M,M) = 0$.
Then the number of pairwise non-isomorphic 
indecomposable direct summands of $M$ is at most $|R^+|$. 
\end{Cor}\vspace{0.2cm}

Let $U_v^-$ be the negative part of the 
quantized enveloping algebra of the Lie algebra corresponding to 
$Q$.
We regard $U_v^-$ as a $\Q(v)$-algebra.
Let ${\mathcal B}$ be the canonical basis and
${\mathcal B}^*$ the dual canonical basis of $U_v^-$,
see \cite{BZ}, \cite{LNT}, \cite{L2} or \cite{R} for definitions.
By \cite[Section 5]{KS}, \cite{L1} the elements of ${\mathcal B}$ (and thus
of ${\mathcal B}^*$) 
correspond to the irreducible components of the nilpotent varieties
$\Lambda({\bf d})$, ${\bf d} \in \N^n$.
Let $b^*(C)$ be the dual canonical basis vector corresponding to
an irreducible component $C$.
We denote the structure constants of $U_v^-$ with respect to the basis
${\mathcal B}^*$ by $\lambda_{C,D}^E$, i.e.
\[
b^*(C)b^*(D) = \sum_E \lambda_{C,D}^E b^*(E).
\]
Following the terminology in \cite{BZ}, two
dual canonical basis vectors
$b^*(C)$ and $b^*(D)$ are called 
{\it multiplicative} if 
\[
b^*(C)b^*(D) = \lambda_{C,D}^E b^*(E)
\] 
for some irreducible component $E$,
and they are {\it quasi-commutative} if
\[
b^*(C)b^*(D) = \lambda b^*(D)b^*(C)
\]
for some $\lambda \in \Q(v)$.
The following conjecture was stated in \cite[Section 1]{BZ}:

\vspace{0.2cm}\begin{Conj}[Berenstein, Zelevinsky]\label{bzconj}
Two dual canonical basis vectors are multiplicative
if and only if they are quasi-commutative.
\end{Conj}\vspace{0.2cm}

One direction of this conjecture was proved by Reineke 
\cite[Corollary 4.5]{R}.
The other direction turned out to be wrong.
Namely,
Leclerc \cite{Le} constructed examples of quasi-commutative 
elements in ${\mathcal B}^*$
which are not multiplicative.
Using preprojective algebras, we give a module theoretic interpretation
of one of his examples. 

Marsh and Reineke \cite{MR} conjectured
that the multiplicative behaviour of
dual canonical basis vectors
should be related to sets of irreducible components with ${\rm Ext}^1$ 
vanishing generically between them.
This was the principal motivation for our work.

The paper is organized as follows: In Section \ref{section2} we review
the main results from \cite{CBSc}.
In Section \ref{section3} we recall known results for the case
that $\Lambda$ is an algebra of finite or tame
representation type.
The proof of Theorem \ref{mainresult} and its corollary can be found in
Section \ref{section4}.
Finally, Section \ref{section5} is devoted to the interpretation of 
Leclerc's example.

\vspace{0.5cm}

{\bf Acknowledgements.} 
The second author is greatful to the IM UNAM, Mexico City, where
most of this work was done. 
We thank Robert Marsh and Markus Reineke 
for helpful and interesting discussions.


\section{Varieties of modules - definitions and known results}
\label{section2}


In this section, we work with arbitrary finite quivers.

\subsection{}
Let
$Q = (Q_0,Q_1)$ be a finite quiver, 
where  $Q_0$ denotes the set of 
vertices and $Q_1$ the set of arrows of $Q$.
Assume that $|Q_0| = n$.
For an arrow $\alpha$ let $s \alpha$ be its starting vertex and
$e \alpha$ its end vertex.
An element ${\bf d} = (d_i)_{i \in Q_0} \in \N^n$ is called
a {\it dimension vector} for $Q$.
A {\it representation} of $Q$ with dimension vector ${\bf d}$ is a matrix
tuple $M = (M_\alpha)_{\alpha \in Q_1}$ with $M_\alpha \in
\M_{d_{s\alpha} \times d_{e\alpha}}(\field)$.
A {\it path} of length $l \geq 1$ in $Q$ is a 
sequence $p = \alpha_1 \cdots \alpha_l$ of arrows in $Q_1$ 
such that $e\alpha_i  = s\alpha_{i+1}$ for $1 \leq i \leq l-1$.
Define $sp = s\alpha_1$ and $ep = e\alpha_l$.
For a representation $M$ and a path $p = \alpha_1 \cdots \alpha_l$
define $M_p = M_{\alpha_1} \cdots M_{\alpha_l}$ which is a matrix
in $\M_{d_{sp} \times d_{ep}}(\field)$.
A {\it relation} for $Q$ is a $\field$-linear combination 
$\sum_{i=1}^t \lambda_i p_i$ of paths $p_i$ of length
at least two such that $sp_i = sp_j$ and $ep_i = ep_j$ for all
$1 \leq i,j \leq t$.
A representation $M$ satisfies such a relation if
$\sum_{i=1}^t \lambda_i M_{p_i} = 0$.
Given a set $\rho$ of relations for $Q$ let $\rep_{(Q,\rho)}({\bf d})$
be the affine variety of representations of $Q$ with dimension
vector ${\bf d}$ which satisfy all relations in $\rho$.

\subsection{}
One can interpret this construction in a module theoretic way.
Namely, let $\field Q$ be the path algebra of $Q$, and let
$A = \field Q/(\rho)$, where $(\rho)$ is the ideal generated by
the elements in $\rho$.
Then $\md_A({\bf d}) = \rep_{(Q,\rho)}({\bf d})$ is the affine 
{\it variety of
$A$-modules} with dimension vector ${\bf d}$.
If $A = \field Q/(\rho)$ is finite-dimensional, then 
$A$ is called a {\it basic algebra}.
In this case, the vertices of $Q$ correspond to the isomorphism
classes of simple $A$-modules, and the entry $d_i$, $i \in Q_0$, 
of ${\bf d}$ is the multiplicity of the simple module corresponding to
$i$ in a composition series of any $M \in \md_A({\bf d})$.
The group $\GL({\bf d}) = \prod_{i \in Q_0} \GL_{d_i}(\field)$ acts
on $\md_A({\bf d})$ by conjugation, i.e. 
\[
g \cdot M = (g_{s\alpha} M_\alpha g_{e\alpha}^{-1})_{\alpha \in Q_1}.
\]
The orbit of $M$ under this action is denoted by $\Or(M)$.
There is a 1-1 correspondence between the set of orbits in $\md_A({\bf d})$ 
and the set of isomorphism classes of $A$-modules with dimension vector 
${\bf d}$.

\subsection{}
Given 
irreducible components $C_i \subseteq \md_A({\bf d}_i)$, $1 \leq i \leq t$,
we consider all $A$-modules with dimension vector ${\bf d} = 
{\bf d}_1 + \cdots + {\bf d}_t$,
which are of the form $M_1 \oplus \cdots \oplus M_t$ with 
the $M_i$ in $C_i$, and we denote by 
$C_1 \oplus \cdots \oplus C_t$ the corresponding
subset of $\md_A({\bf d})$. 
This is the image of the map
\[
\GL({\bf d}) \times C_1 \times \cdots \times C_t
\longrightarrow \md_A({\bf d}) 
\]
\[
(g,M_1,\cdots,M_t) \mapsto g \cdot \left(
\bigoplus_{i=1}^t M_i \right).
\]
We call $C_1 \oplus \cdots \oplus C_t$ the {\it direct sum} of the 
components $C_i$.
It follows that the
closure $\overline{C_1 \oplus \cdots \oplus C_t}$ is irreducible.
For an irreducible component $C$ define $C^n = \bigoplus_{i=1}^n C$.
We call $C$ {\it indecomposable} if $C$ contains a dense subset of
indecomposable $A$-modules.
The following result from \cite{CBSc} is
an analogue of the Krull-Remak-Schmidt
Theorem.

\vspace{0.2cm}\begin{Thm}\label{thm1}
If $C \subseteq \md_A({\bf d})$ is an irreducible component,
then 
\[
C = \overline{C_1 \oplus \cdots \oplus C_t}
\]
for some indecomposable 
irreducible components $C_i \subseteq \md_A({\bf d}_i)$,
$1 \leq i \leq t$, and
$C_1, \cdots, C_t$
are uniquely determined by this, up to reordering.
The above direct sum is called the canonical decomposition of $C$.
\end{Thm}\vspace{0.2cm}

However, the closure of a direct sum of irreducible
components is not in general an irreducible component.
The next result is also proved in \cite{CBSc}.

\vspace{0.2cm}\begin{Thm}\label{thm2}
If $C_i \subseteq \md_A({\bf d}_i)$, $1 \leq i \leq t$, are irreducible
components, and ${\bf d} = {\bf d}_1 + \cdots + {\bf d}_t$,
then $\overline{C_1 \oplus \cdots \oplus C_t}$
is an irreducible component of $\md_A({\bf d})$
if and only if ${\rm ext}_A^1(C_i,C_j) = 0$ for all $i \not= j$.
\end{Thm}\vspace{0.2cm}

Instead of taking direct sums of the modules in two irreducible
components, one can take extensions. 
Let ${\bf d} = {\bf d}_1 + {\bf d}_2$ be dimension vectors, let
$G = \GL({\bf d}_1) \times \GL({\bf d}_2)$,
and let $S$ be a $G$-stable subset of 
$\md_A({\bf d}_1) \times \md_A({\bf d}_2)$.
We denote by ${\mathcal E}(S)$
the $\GL({\bf d})$-stable subset of $\md_A({\bf d})$
corresponding to all modules $M$ which belong to
a short exact sequence
\[
0 \longrightarrow M_2 \longrightarrow M \longrightarrow M_1 
\longrightarrow 0
\]
with $(M_1,M_2) \in S$, see \cite{CBSc} for more details.

For an irreducible component $C \subseteq \md_A({\bf d})$ let
\[
\mu_g(C) = {\rm dim}\, C - \max \{ {\rm dim}\, \Or(M) \mid M \in C \}
\]
be the {\it generic number of parameters of} $C$.
Thus $\mu_g(C) = 0$ if and only if $C$ contains a dense orbit $\Or(M)$.
For example, if $P$ is a projective $A$-module, then 
$\Ext_A^1(P,P) = 0$.
This implies that the closure of the orbit $\Or(P)$ is an irreducible
component, and we get $\mu_g(\overline{\Or(P)}) = 0$.
Also, if $C = \overline{C_1 \oplus \cdots \oplus C_t}$ with
${\rm ext}^1_A(C_i,C_j) = 0$ for all $i \not= j$, then
\[
\mu_g(C) = \sum_{i=1}^t \mu_g(C_i).
\]


\section{The finite and tame cases}
\label{section3}


As in the introduction let $Q$ be a Dynkin quiver.
Then $\Lambda$ is of finite representation 
type if and only if
$Q$ is of type $\A_i$, $i \leq 4$.
In this case, 
if
$\{ C_i \subseteq \Lambda({\bf d}_i) \mid i \in I \}$ is a maximal set of 
pairwise different indecomposable irreducible
components such that ${\rm ext}_{\Lambda}^1(C_i,C_j) = 0$ for all $i,j$,
then $|I| = |R^+|$.
This follows from \cite{BZ} for $i \leq 3$, and the case $i = 4$ was
done by Marsh and Reineke.

Recall that for a tame algebra $A$ one has
$\mu_g(C) \leq 1$ for any indecomposable irreducible
component $C \subseteq \md_A({\bf d})$.
It is known that $\Lambda$ is of tame representation type
if and only if $Q$ is of type $\A_5$ or $\D_4$.
In this case, a complete classification of the indecomposable
irreducible components, and a necessary and sufficient condition
for ${\rm ext}^1_{\Lambda}(C,D) = 0$ for any two irreducible components
$C$ and $D$ was obtained in \cite{GSc}.
In particular, this implies the following:

\vspace{0.2cm}\begin{Thm}
Assume that $Q$ is of type $\A_5$ or $\D_4$.
Then the following hold:
\begin{itemize}

\item[(1)] For any irreducible component $C \subseteq \Lambda({\bf d})$
we have ${\rm ext}_{\Lambda}^1(C,C) = 0$;

\item[(2)] If $C \subseteq \Lambda({\bf d})$ is an indecomposable
irreducible component, then we have $\mu_g(C) = 0$ or
$\mu_g(C) = 1$.
For suitable ${\bf d}$ there exists an indecomposable irreducible
component $C \subseteq \Lambda({\bf d})$ with
$\mu_g(C) = 1$;

\item[(3)] 
Let $\{ C_i \subseteq \Lambda({\bf d}_i) 
\mid i \in I \}$ 
be a maximal set of pairwise different indecomposable irreducible 
components such that
${\rm ext}_{\Lambda}^1(C_i,C_j) = 0$ for all $i,j$.
Then there is at most one $C_i$ with $\mu_g(C_i) = 1$.
In this case, we have
$|I| = |R^+| - 1$, 
and we get
$|I| =  |R^+|$, otherwise.

\end{itemize}
\end{Thm}\vspace{0.2cm}

This leads us to the following conjecture for arbitrary Dynkin quivers of 
type $\A_n$, $\D_n$ or $\E_{6,7,8}$:

\vspace{0.2cm}\begin{Conj}\label{Conj1}
If
$\{ C_i \subseteq \Lambda({\bf d}_i) \mid i \in I \}$ is a maximal set of 
pairwise different indecomposable irreducible
components such that ${\rm ext}_{\Lambda}^1(C_i,C_j) = 0$ for all $i,j$,
then
\[
|I| = |R^+| - \sum_{i \in I} \mu_g(C_i).
\] 
\end{Conj}\vspace{0.2cm}

In all remaining cases the algebra $\Lambda$ is
of wild representation type.
So one should expect 
irreducible components $C$ with ${\rm ext}_{\Lambda}^1(C,C) \not= 0$.
Thus, maybe one should study sets 
$\{ C_i \subseteq \Lambda({\bf d}_i) \mid i \in I \}$ of
irreducible components with the weaker condition 
${\rm ext}_{\Lambda}^1(C_i,C_j) = 0$
for all $i \not= j$.
However, we do not know how to generalize Theorem \ref{mainresult} to
this case.


\section{Proof of Theorem \ref{mainresult}}
\label{section4}


As before let $Q$ be a Dynkin quiver, and 
let $R^+ = \{ {\bf a}_i \mid 1 \leq i \leq N \}$
be the set of positive roots of $Q$.
By Gabriel's Theorem there is a 1-1 correspondence between the
isomorphism classes of indecomposable $\field Q$-modules and
the elements in $R^+$.
This correspondence associates to a root ${\bf a}_i$ the
isomorphism class $[M({\bf a}_i)]$ of an
indecomposable $\field Q$-module $M({\bf a}_i)$ with dimension
vector ${\bf a}_i$.
By the Theorem of Krull-Remak-Schmidt each $\field Q$-module is 
isomorphic to a (up to reordering)
unique direct sum of the indecomposable modules $M({\bf a}_i)$.
The maps
\[
\alpha = (\alpha_1, \cdots, \alpha_N) \mapsto
\Or(M_\alpha) \mapsto C_\alpha \mapsto b^*(C_\alpha)
\]
define 1-1 correspondences between $\N^N$, the set of orbits
$\Or(M) \subseteq \md_{\field Q}({\bf d})$, ${\bf d} \in \N^n$,
the set of irreducible components of 
$\Lambda({\bf d})$, ${\bf d} \in \N^n$,
and the set ${\mathcal B}^*$ of dual canonical basis vectors,
where
\begin{eqnarray*}
M_\alpha &=& \bigoplus_{i=1}^N M({\bf a}_i)^{\alpha_i},\\
C_\alpha &=&  \overline{\pi^{-1}(\Or(M_\alpha))}
\end{eqnarray*} 
with
\[
\pi: \Lambda({\bf d}) \to \md_{\field Q}({\bf d})
\]
the canonical projection map.

Let $\alpha, \beta \in \N^N$.
By Theorem \ref{thm2} the closure
$\overline{C_\alpha \oplus C_\beta}$ is an irreducible
component if and only if ${\rm ext}^1_{\Lambda}(C_\alpha,C_\beta)
= {\rm ext}^1_{\Lambda}(C_\beta,C_\alpha) = 0$.
In this case, we have $\overline{C_\alpha \oplus C_\beta}
= C_{\alpha + \beta}$.

Let $\alpha_i = (\alpha_{1i}, \cdots, \alpha_{Ni})$, 
$1 \leq i \leq N+1$, be non-zero pairwise different elements in $\N^N$
such that $C_{\alpha_i}$ is an indecomposable irreducible component
for all $i$.
To get a contradiction, we assume that
${\rm ext}^1_{\Lambda}(C_{\alpha_i},C_{\alpha_j}) = 0$ for all
$1 \leq i,j \leq N+1$.
For all ${\bf m} = (m_1, \cdots, m_{N+1}) \in \N^{N+1}$ define
\[
C_{\bf m} = C_{\sum_{i=1}^{N+1} m_i \alpha_i}.
\]
We get
\[
C_{\bf m} = 
\overline{C_{\alpha_1}^{m_1} \oplus \cdots \oplus 
C_{\alpha_{N+1}}^{m_{N+1}}}.
\]
Since the $C_{\alpha_i}$ are indecomposable, the above is the canonical
decomposition of the irreducible component 
$C_{\bf m}$.

We claim that there exist some elements 
${\bf m} = (m_1, \cdots, m_{N+1})
\not= {\bf l} = (l_1, \cdots, l_{N+1})$ in $\N^{N+1}$
such that
\[
\sum_{i=1}^{N+1} m_i \alpha_i = \sum_{i=1}^{N+1} l_i \alpha_i.
\]
This implies $C_{\bf m} = C_{\bf l}$.
Thus, we get a contradiction to the unicity 
of the canonical decomposition of irreducible components, 
see Theorem \ref{thm1}.

Let
$\Delta = (\alpha_{ij})$ be the $N \times (N+1)$-matrix where the 
$j$th column
is just the vector $\alpha_j$.
Thus we have to find some ${\bf d} \in \N^N$ and some
${\bf m} \not= {\bf l}$ in $\N^{N+1}$ such that
\[
\Delta{\bf m} = \Delta{\bf l} = {\bf d}.
\]
Since all entries in $\Delta$ are in $\N$, this would imply 
that ${\bf m}$, ${\bf l}$ and ${\bf d}$ are all non-zero.
Furthermore, there must be a non-zero element 
${\bf z} = (z_1, \cdots, z_{N+1}) 
\in \Z^{N+1}$ such that $\Delta{\bf z} = 0$.

First, we consider the case ${\bf z} \in \N^{N+1}$.
Let ${\bf m} = (1, \cdots, 1) \in \N^{N+1}$, ${\bf d} = \Delta{\bf m}$ and
${\bf l} = {\bf m} + {\bf z}$.
Obviously, ${\bf d} \in \N^N$.
We get $\Delta{\bf m} = \Delta{\bf l} = {\bf d}$ with 
${\bf m} \not= {\bf l}$ in $\N^{N+1}$.

Next, assume that ${\bf z} \notin \N^{N+1}$.
Define $\lambda = - \min\{ z_i \mid 1 \leq i \leq N+1 \}$.
Let ${\bf m} = (\lambda, \cdots, \lambda) \in \N^{N+1}$,
${\bf d} = \Delta{\bf m}$ and
${\bf l} = {\bf m} + {\bf z}$.
Again, we get  ${\bf d} \in \N^N$ and
$\Delta{\bf m} = \Delta{\bf l} = {\bf d}$ with
${\bf m} \not= {\bf l}$ in $\N^{N+1}$.
This finishes the proof of Theorem \ref{mainresult}.

Corollary \ref{maincor} follows immediately from
the fact that an orbit ${\mathcal O}(N) \subseteq \md_A({\bf d})$
of an $A$-module $N$ is open provided ${\rm Ext}_A^1(N,N) = 0$.
Clearly, ${\mathcal O}(N)$ is open if and only if the closure
$\overline{{\mathcal O}(N)}$ is an irreducible component.
Then we use Theorems \ref{mainresult} and \ref{thm2}.


\section{Interpretation of Leclerc's example}\label{section5}


In the following, we use the notation introduced at the beginning
of Section \ref{section4}.

Reineke proved in \cite[Lemma 4.6]{R} that the multiplicativity of 
$b^*(C_\alpha)$ and $b^*(C_\beta)$ implies that 
\[
b^*(C_\alpha)b^*(C_\beta) = v^mb^*(C_{\alpha+\beta})
\]
for some $m \in \Z$.
He also showed that $\lambda_{C,D}^E \not= 0$ if and only if
$\lambda_{D,C}^E \not= 0$.
This follows from \cite[Proposition 4.4]{R}.
Thus one direction of Conjecture \ref{bzconj} holds,
namely if two dual canonical basis vectors are multiplicative, then
they are quasi-commutative.
The following related problem should be of interest:

\vspace{0.2cm}\begin{Prob}
Describe the elements $\alpha, \beta \in \N^N$ such that
\[
C_{\alpha+\beta} = \overline{C_\alpha \oplus C_\beta}.
\]
\end{Prob}\vspace{0.2cm}

As mentioned in the introduction, 
Leclerc recently constructed in \cite{Le} counterexamples for the
other direction of the Berenstein-Zelevinsky Conjecture.
We give a module theoretic interpretation of one of his examples:

Let $Q$ be the quiver of type $\A_5$
with arrows $a_i: i+1 \to i$, $1 \leq i \leq 4$.
Thus $\Lambda$ is given by the quiver 
\[
\xymatrix{
1\ar@<1ex>[r]^{\bar{a}_1}&\ar@<1ex>[l]^{a_1} 2 \ar@<1ex>[r]^{\bar{a}_2}&
\ar@<1ex>[l]^{a_2} 3\ar@<1ex>[r]^{\bar{a}_3}&\ar@<1ex>[l]^{a_3} 4 
\ar@<1ex>[r]^{\bar{a}_4}&
\ar@<1ex>[l]^{a_4} 5
}
\]
and the following set of relations
\[
\{  \bar{a}_1a_1, a_1\bar{a}_1-\bar{a}_2a_2,a_2\bar{a}_2-\bar{a}_3a_3,
a_3\bar{a}_3-\bar{a}_4a_4,a_4\bar{a}_4 \}.
\]
Now $R^+$ contains exactly 15 elements, namely for each
$1 \leq i \leq j \leq 5$ there is a positive root 
$[i,j] = (d_l)_{1 \leq l \leq 5}$ with $d_l = 1$ for
$i \leq l \leq j$, and $d_l = 0$, else.
We identify $\N R^+$ with 
$\N^{15}$ by fixing a linear ordering on
$R^+$, namely let
\begin{eqnarray*}
& & [1,1] < [1,2] < [1,3]<[1,4]<[1,5]<[2,2]<[2,3]<[2,4]<[2,5]\\
& & <[3,3]<[3,4]<[3,5]<[4,4]<[4,5]<[5,5].
\end{eqnarray*}
Define
\begin{eqnarray*} 
\alpha &=& [1,2] + [2,4] + [3,3] + [4,5],\\
\beta &=& [1,2] + [1,4] +  [2,3] + [2,5] + [3,4] + 
[4,5].
\end{eqnarray*}
Thus, regarded as elements in $\N^{15}$ we have
\begin{eqnarray*}
\alpha &=& (0,1,0,0,0,0,0,1,0,1,0,0,0,1,0),\\
\beta &=&  (0,1,0,1,0,0,1,0,1,0,1,0,0,1,0).
\end{eqnarray*}
In \cite{Le}
Leclerc showed that
\[
b^*(C_\alpha)^2 = v^{-2} \left( b^*(C_{\alpha+\alpha})
+ b^*(C_\beta) \right).
\]
This is obviously a counterexample to the Berenstein-Zelevinsky Conjecture.
Now define  
\begin{eqnarray*} 
\beta_1 &=& [1,2] + [2,3] + [3,4] + [4,5],\\ 
\beta_2 &=& [1,4] + [2,5].
\end{eqnarray*}
Thus we have $\beta = \beta_1 + \beta_2$.

\vspace{0.2cm}\begin{Prop}
Let $\alpha, \beta, \beta_1,\beta_2$ be as above.
Then the following hold:
\begin{itemize}

\item[(1)]
The irreducible components $C_\alpha$, $C_{\beta_1}$ and
$C_{\beta_2}$ are indecomposable with
$\mu_g(C_\alpha) = 1$ and $\mu_g(C_{\beta_1}) = \mu_g(C_{\beta_2}) = 0$;

\item[(2)]
We have
$C_{\alpha+\alpha} = \overline{C_\alpha \oplus C_\alpha}$ and
$C_\beta = \overline{C_{\beta_1} \oplus C_{\beta_2}}$.
Thus 
\[
b^*(C_\alpha)^2 = 
v^{-2} \left( b^*(\overline{C_{\alpha} \oplus 
C_{\alpha}})
+ b^*(\overline{C_{\beta_1} \oplus C_{\beta_2}}) \right).
\]
\end{itemize}
\end{Prop}

\begin{proof}
For $\lambda \in \field \setminus \{ 0,1 \}$ let $M_\lambda$ be 
the $8$-dimensional $\Lambda$-module 
where the arrows of $\Lambda$ operate on a basis $\{ 1, \cdots, 8 \}$
as in the following picture:

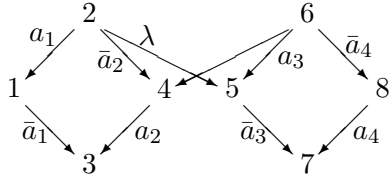
\begin{center}
\unitlength0.5cm
\begin{picture}(11,6)

\put(3,5){$2$}
\put(2.8,4.8){\vector(-1,-1){1.3}}\put(1.6,4.5){$a_1$}
\put(3.5,4.8){\vector(1,-1){1.3}}\put(3.4,3.8){$\bar{a}_2$}
\put(3.6,4.8){\vector(2,-1){3}}\put(4.5,4.4){$\lambda$}
\put(1,3){$1$}
\put(1.5,2.8){\vector(1,-1){1.3}}\put(1.4,1.9){$\bar{a}_1$}
\put(3,1){$3$}
\put(5,3){$4$}
\put(4.8,2.8){\vector(-1,-1){1.3}}\put(4.4,1.9){$a_2$}

\put(8.8,5){$6$}
\put(8.6,4.8){\vector(-1,-1){1.3}}\put(8.2,3.9){$a_3$}
\put(9.3,4.8){\vector(1,-1){1.3}}\put(10,4.2){$\bar{a}_4$}
\put(8.5,4.8){\vector(-2,-1){3}}
\put(6.8,3){$5$}
\put(7.3,2.8){\vector(1,-1){1.3}}\put(7.2,1.9){$\bar{a}_3$}
\put(8.8,1){$7$}
\put(10.8,3){$8$}
\put(10.6,2.8){\vector(-1,-1){1.3}}\put(10.2,1.9){$a_4$}

\end{picture}
\end{center}

Thus, for example $2 \cdot a_1 = 1$,
$2 \cdot \bar{a}_2 = 4 + \lambda 5$, $6 \cdot a_3 = 4 + 5$,
$1 \cdot \bar{a}_1 = 3$, etc.
Note that $M_\lambda$ lies in $C_\alpha$.

The modules $M_\lambda$ are indecomposable and
$
{\rm dim}\, {\rm End}_{\Lambda}(M_\lambda) = 3.
$
From a well-known general fact we know that
each irreducible component of $\Lambda(1,2,2,2,1)$
has dimension $2+4+4+2 = 12$, see for example \cite[Section 12]{L2}.
The group $\GL(1,2,2,2,1)$ acts as described in Section \ref{section2}
on $\Lambda(1,2,2,2,1)$ and has dimension $14$.
Thus we get
$
{\rm dim}\, {\mathcal O}(M_\lambda) = 14-3 = 11.
$
One checks easily that $M_\lambda$ and $M_\mu$ are isomorphic if and
only if $\lambda = \mu$.
This implies
\[
{\rm dim}\, \overline{ \{ {\mathcal O}(M_\lambda) \mid
\lambda \in \field \setminus \{ 0,1 \} \} } = 11+1 = 12.
\]
We get
\[
C_\alpha =  \overline{\{ {\mathcal O}(M_\lambda) \mid
\lambda \in \field \setminus \{ 0,1 \} \}}.
\]
Thus $C_\alpha$ is an indecomposable irreducible component
with $\mu_g(C_\alpha) = 1$.

Next, let $P_2$ and $P_4$ be the indecomposable projective 
$\Lambda$-modules corresponding to the vertices $2$ and $4$, respectively.
These modules are both $8$-dimensional and 
look as in the following picture:

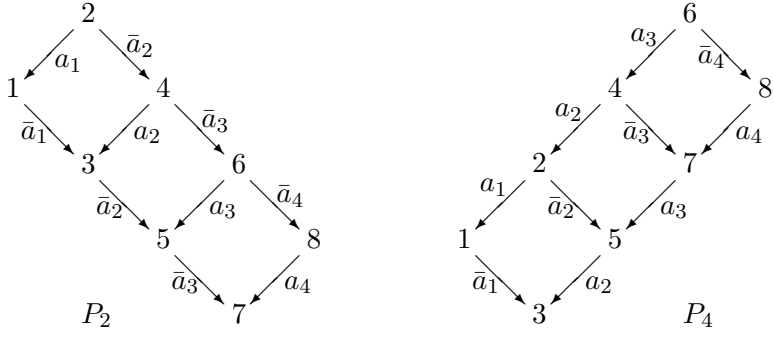
\begin{center}
\unitlength0.5cm
\begin{picture}(22,10)

\put(3,9){$2$}
\put(2.8,8.8){\vector(-1,-1){1.3}}\put(2.3,7.8){$a_1$}
\put(3.5,8.8){\vector(1,-1){1.3}}\put(4.2,8.2){$\bar{a}_2$}
\put(1,7){$1$}
\put(1.5,6.8){\vector(1,-1){1.3}}\put(1.4,5.9){$\bar{a}_1$}
\put(3,5){$3$}
\put(3.5,4.8){\vector(1,-1){1.3}}\put(3.4,3.9){$\bar{a}_2$}
\put(3,1){$P_2$}
\put(5,7){$4$}
\put(4.8,6.8){\vector(-1,-1){1.3}}\put(4.4,5.9){$a_2$}
\put(5.5,6.8){\vector(1,-1){1.3}}\put(6.2,6.2){$\bar{a}_3$}
\put(5,3){$5$}
\put(5.5,2.8){\vector(1,-1){1.3}}\put(5.4,1.9){$\bar{a}_3$}
\put(7,5){$6$}
\put(6.8,4.8){\vector(-1,-1){1.3}}\put(6.4,3.9){$a_3$}
\put(7.5,4.8){\vector(1,-1){1.3}}\put(8.2,4.2){$\bar{a}_4$}
\put(7,1){$7$}
\put(9,3){$8$}
\put(8.8,2.8){\vector(-1,-1){1.3}}\put(8.4,1.9){$a_4$}

\put(13,3){$1$}
\put(13.5,2.8){\vector(1,-1){1.3}}\put(13.4,1.9){$\bar{a}_1$}
\put(15,5){$2$}
\put(14.8,4.8){\vector(-1,-1){1.3}}\put(13.6,4.5){$a_1$}
\put(15.5,4.8){\vector(1,-1){1.3}}\put(15.4,3.9){$\bar{a}_2$}
\put(17,7){$4$}
\put(16.8,6.8){\vector(-1,-1){1.3}}\put(15.6,6.5){$a_2$}
\put(17.5,6.8){\vector(1,-1){1.3}}\put(17.4,5.9){$\bar{a}_3$}
\put(19,9){$6$}
\put(18.8,8.8){\vector(-1,-1){1.3}}\put(17.6,8.5){$a_3$}
\put(19.5,8.8){\vector(1,-1){1.3}}\put(19.4,7.9){$\bar{a}_4$}
\put(15,1){$3$}
\put(17,3){$5$}
\put(16.8,2.8){\vector(-1,-1){1.3}}\put(16.4,1.9){$a_2$}
\put(19,5){$7$}
\put(18.8,4.8){\vector(-1,-1){1.3}}\put(18.4,3.9){$a_3$}
\put(19,1){$P_4$}
\put(21,7){$8$}
\put(20.8,6.8){\vector(-1,-1){1.3}}\put(20.4,5.9){$a_4$}

\end{picture}
\end{center}

We have ${\rm Ext}_{\Lambda}^1(P_i,P_j) = 0$ for
all $i,j \in \{ 2,4 \}$.
This follows directly from the projectivity of both modules.
From this and the above pictures we get
\begin{eqnarray*}
C_{\beta_1} &=& \overline{{\mathcal O}(P_2)},\\ 
C_{\beta_2} &=& \overline{{\mathcal O}(P_4)},\\
C_\beta &=& \overline{C_{\beta_1} \oplus C_{\beta_2}}.
\end{eqnarray*}
In particular, $C_{\beta_1}$ and $C_{\beta_2}$ are indecomposable
irreducible components with $\mu_g(C_{\beta_1}) = \mu_g(C_{\beta_2}) = 0$.
This finishes the proof.
\end{proof}

For irreducible components $C \subseteq \Lambda({\bf d})$ and 
$D \subseteq \Lambda({\bf e})$ define
\[
{\mathcal V}(C,D) = 
\bigcap_{U \subseteq C, V \subseteq D} 
\left\{ E \subseteq \Lambda({\bf d}+{\bf e}) \text{ irred. comp.} \mid 
E \subseteq 
\overline{{\mathcal E}(U \times V)} \right\},
\]
where $U$ (resp. $V$) runs through all non-empty $\GL({\bf d})$-stable 
(resp. $\GL({\bf e})$-stable) open subsets of
$C$ (resp. $D$),
see Section \ref{section2} for the definition of ${\mathcal E}(U \times V)$.

Using the previous proposition, and some well-known results
on the representation theory of the algebra $\Lambda$,
see \cite{GSc} and \cite{Ri},
one can show that
\[
{\mathcal V}(C_\alpha,C_\alpha) = \{ C_{\alpha+\alpha},
C_\beta \}.
\]
Note that ${\rm ext}_{\Lambda}^1(C_\alpha,C_\alpha) = 0$,
since ${\rm Ext}_{\Lambda}^1(M_\lambda,M_\mu) = 0$ for
all $\lambda \not= \mu$. 
But one can show that 
${\rm dim}\, {\rm Ext}_{\Lambda}^1(M_\lambda,M_\lambda) = 2$, 
see \cite[Section 6]{GSc}.
For any $M_\lambda$ there is a short exact sequence
\[
0 \to M_\lambda \to M_\lambda(2) \to M_\lambda \to 0,
\]
where $M_\lambda(2)$ is the module of quasi-length two in the same 
Auslander-Reiten component as $M_\lambda$ (it is known that
$M_\lambda$ lies in a homogeneous tube).
Additionally to this `natural' self-extension,  
there exists a short exact sequence
\[
0 \to M_\lambda \to P_2 \oplus P_4 \to M_\lambda \to 0.
\]
Motivated by our above analysis, 
one might ask whether the following is true:

\vspace{0.2cm}\begin{Qu}
Are the following two statements equivalent?
\begin{itemize}

\item[(1)]
$\lambda_{C,D}^E \not= 0$;

\item[(2)] 
$E \in {\mathcal V}(C,D) \cup {\mathcal V}(D,C)$.

\end{itemize}
\end{Qu}\vspace{0.2cm}

A positive answer to the above question would imply the following:
\begin{itemize}

\item 
If an irreducible component $C$ contains an open orbit, then
\[
b^*(C)^2 = v^m b^*(\overline{C \oplus C})
\]
for some $m \in \Z$.
Here we use that ${\rm Ext}^1_{\Lambda}(M,M) = 0$ if and only if
${\mathcal O}(M)$ is an open orbit.
This is a special feature of preprojective algebras.

\item
If irreducible components $C$ and $D$ contain non-empty stable open
subsets $U \subseteq C$ and $V \subseteq D$ such that
${\rm Ext}_{\Lambda}^1(M,N) = 0$ for all $M \in U, N \in V$, then
$b^*(C)$ and $b^*(D)$ are multiplicative.
Here we use \cite[Theorem 1.3]{CBSc}.

\item
If an irreducible component $C$ contains an open orbit, and
if $D$ is an irreducible component such that
${\rm ext}_{\Lambda}^1(C,D) = 0$, then
$b^*(C)$ and $b^*(D)$ are multiplicative.
Again this uses \cite[Theorem 1.3]{CBSc}.

\end{itemize}


\end{document}